\documentclass{article}

\usepackage{PRIMEarxiv}

\usepackage[utf8]{inputenc} % allow utf-8 input
\usepackage[T1]{fontenc}    % use 8-bit T1 fonts
\usepackage{hyperref}       % hyperlinks
\usepackage{url}            % simple URL typesetting
\usepackage{booktabs}       % professional-quality tables
\usepackage{amsfonts}       % blackboard math symbols
\usepackage{microtype}      % microtypography
\usepackage{lipsum}
\usepackage{fancyhdr}       % header
\usepackage{graphicx}       % graphics

\usepackage{amsfonts}
\usepackage{amssymb}
\usepackage{amsthm}
\graphicspath{{media/}}     % organize your images and other figures under media/ folder

\newtheorem{theorem}{Theorem}
\newtheorem{definition}{Definition}

\newtheorem{remark}{Remark}
\newtheorem{lemma}{Lemma}
\newtheorem{Example}{Example}
\newtheorem{corollary}{Corollary}

\title{Fuzzy Volterra Integral Equations with Piecewise Continuous Kernels: Theory and Numerical Solution
}

\author{
  Samad Noeiaghdam, \, Aliona Dreglea \\
  Industrial Mathematics Lab\\
  Baikal School of BRICS\\
   Irkutsk National Research Technical University \\
  Irkutsk\\
  \texttt{\{snoei, adreglea\}@istu.edu} \\
   \And
  Denis N. Sidorov \\
  Department of Applied Mathemtics \\
  Melentiev Energy Systems Institute\\ 
  Siberian Branch of Russian Academy of Sciences\\
  Irkutsk\\
  \texttt{dsidorov@isem.irk.ru} \\
}

\begin{document}
\maketitle

\begin{abstract}
This study aims to discuss the existence and uniqueness of solution of fuzzy Volterra integral equations  with piecewise continuous kernels.  These types of problems are often encountered in balancing issues for systems with hereditary dynamics, such as electric load leveling. The method of successive approximations is applied and the main theorems are proved based on the method. Some examples are discussed and the results are presented for different values of $\mu$ by plotting several graphs.
\end{abstract}

% keywords can be removed

\section{Introduction}

Fuzzy integral equations (FIEs) are among applicable and important problems of engineering and basic sciences.  Bede and Gal \cite{m8}, Friedman and Ma \cite{m18} and Goetschel and Voxman \cite{m22} have some studies on theory of FEIs. Ziari and Abbasbandy solved nonlinear FEIs using fuzzy quadrature rules \cite{1}. The Reproducing Kernel Hilbert space method has been applied by Javan et al. in \cite{2}, the radial basis functions has discussed by Asari et al in \cite{3}. Amirfakhrian et al. used the fuzzy interpolation techniques for solving FIEs in \cite{4}. Also many other techniques for solving FIEs can be found in \cite{5}. In \cite{13} the well-known sinc-collcation method in both DE and SE precisions were used for solving fuzzy Fredholm integral equations.  In \cite{16} combining of the homotopy analysis method and  Laplace transformations were applied to study the FEIs of the Abel type. In \cite{17,18} the CESTAC method and the CADNA library were employed to identify the optimal results of the homotopy analysis method for solving FIEs.

Volterra integral equation with piecewise continuous kernel is known and applicable problem which can be 
employed in various balance problems including electric loading problem. Sidorov et al. in \cite{25,26} studied the generalized solution of Volterra integral equations. Solvability of this problem has been illustrated by Sidorov in \cite{28,29} and Muftahov and Sidorov in \cite{27}. The successive approximation  method was used to find the solution of Volterra integral equations in \cite{31}. The numerical solution of this problem can be found in \cite{32}.  Also some numerical and semi-analytical methods can be found for solving Volterra integral equations with piecewise kernel such as the spline collocation method \cite{19}, Lagrange-collocation method \cite{20}, Adomian decomposition method \cite{22}, homotopy perturbation method \cite{23}, the collocation method with Taylor polynomials \cite{24} and other \cite{21}. For more details on the theory of Volterra integral equations 
with piecewise continuous kernels readers may refer to  monograph \cite{29a}. Such equations naturally 
generalizes the non-classic Volterra equations studied in  monograph \cite{29aa}.

This study  deals with the novel class of  fuzzy Voltera integral equation  (FVIE) with piecewise continuous kernel
\begin{equation}\label{4}
Z(v) = Y(v) \oplus (\mathcal{FR}) \sum_{t=1}^{m'} \int_{\theta_{t-1}(v)}^{\theta_{t}(v)} K_t(r,v) \odot G (Z(r)) dr,~~~z_1 \leq s,v \leq T \leq z_2,
\end{equation}
where
$$
z_1=: \theta_0(v) < \theta_1(v) < ... < \theta_{m'-1}(v) < \theta_{m'}(v):=v,~~~ z_1 \leq v \leq T \leq z_2
$$
and the kernel $K_t(r,v)$ is a crisp  and positive function over the square  $z_1 \leq s,v \leq T \leq z_2$, $Z(v)$ shows a fuzzy real valued function and $G: \mathbb{R}_\digamma \rightarrow \mathbb{R}_\digamma$ is continuous. Also $K_t(r,v)$ is a piecewise kernel along continuous curves $\theta_t(v), t=1,2,...,m'$, therefore $K_1(r,v), K_2(r,v), ..., K_{m'}(r,v)$ are uniformly continuous with respect to $t$ and there exist $M_t > 0$ such that $M_t = \max_{z_1 \leq r,v \leq z_2}  | K_t(r,v)  |$.  We applied the successive approximations for solving problem (\ref{4}). The existence of solution theorem is also discussed. Also the main theorem is proved below to show the error estimation of the problem. Solving some examples in both linear and nonlinear and plotting error graphs and also graphs of fuzzy approximate solutions, the ability and efficiency of the method are shown.

This paper is organized as follows. Section 2 provdes the preliminaries of fuzzy mathematics. Section 3 is the main idea of this study. Also in this section  the main existence of solution theorem is illustrated. Section 4 shows the error estimation of the successive approximation method for solving problem (\ref{4}). Section 5 provides the linear and nonlinear  examples. Using some graphs we show the accuracy of the method. Section 5 is the conclusion.

%%%======================================================================================================================
%%%======================================================================================================================
%%%======================================================================================================================
\section{Preliminaries}
We have reported the  main definitions and theorems of fuzzy mathematics \cite{m8,m13,m18,m22,m23,m25}.

\begin{definition} \label{def1}
Based on the following properties a fuzzy number $p : \mathbb{R} \rightarrow [0,1]$ can be defined as a function: 
\begin{enumerate}
  \item $p$ is normal which is $\exists x_0 \in \mathbb{R}; p(x_0)=1$,
  \item $p$ is fuzzy convex set $p(\gamma x +(1-\gamma )y) \geq \min \{ p(x), p(y) \}, \forall x,y \in \mathbb{R}, \gamma \in [0,1] $.
    \item $p$ is upper semi-continuous on $\mathbb{R}$,
    \item $\{ x \in \mathbb{R}: p(x) >0 \}$ is a compact set.
\end{enumerate}
\end{definition}
$\mathbb{R}_\digamma$ shows all fuzzy numbers sets.

\begin{definition} \label{def2}

$(\underline{p}(\mu), \overline{p}(\mu)), 0 \leq \mu \leq 1$ is the parametric form of an arbitrary fuzzy number satisfying the following conditions:
\begin{enumerate}
  \item $\underline{p}(\mu)$ is a bounded left continuous non-decreasing function over $[0, 1]$,
  \item $\overline{p}(\mu)$ is a bounded left continuous non-increasing function over $[0, 1]$,
  \item $\underline{p}(\mu) \leq \overline{p}(\mu), 0 \leq \mu \leq 1$
\end{enumerate}
We show the  scalar multiplication and addition of fuzzy numbers as:
\begin{enumerate}
  \item $(p\oplus p_1)(\mu) = (\underline{p}(\mu)+\underline{p_1}(\mu), \overline{p}(\mu)+\overline{p_1}(\mu))$,
  \item  $ (\gamma \odot p ) (\mu) = \left\{  \begin{array}{ll}
                                               (\gamma \underline{p}(\mu), \gamma \overline{p}(\mu)) & \gamma \geq 0, \\
                                              (\gamma \overline{p}(\mu), \gamma \underline{p}(\mu)) & \gamma < 0.
                                             \end{array} \right.
  $
\end{enumerate}
\end{definition}

\begin{definition} \label{def4}
Let $p=(\underline{p}(\mu), \overline{p}(\mu)), p_1=(\underline{p_1}(\mu), \overline{p_1}(\mu))$ be two fuzzy numbers then the distance can be defined as
$$\mathcal{D}(p,p_1) = \sup_{\mu \in [0,1]} \max \{ |  \underline{p}(\mu) - \underline{p}(\mu) |,  | \overline{p}(\mu) - \overline{p_1}(\mu) |  \}. $$

We have the following properties for distance $\mathcal{D}$.
\end{definition}

\begin{theorem} \label{th1}
\begin{enumerate}
  \item $(\mathbb{R}_\digamma, \mathcal{D})$ is a complete metric space,
  \item $\mathcal{D}(p\oplus p_2, p_1\oplus p_2) = \mathcal{D}(p,p_1)  \forall p, p_1, p_2 \in \mathbb{R}_\digamma$,
  \item $\mathcal{D}(k \odot p, k \odot p_1) = | k| \mathcal{D}(p,p_1), \forall p, p_1 \in \mathbb{R}_\digamma \forall k \in \mathbb{R}$,
  \item $\mathcal{D} (p \oplus p_1, p_2 \oplus p_3) \leq \mathcal{D} (p,p_2) + \mathcal{D} (p_1,p_3) \forall p, p_1, p_2, p_3 \in \mathbb{R}_\digamma$.
\end{enumerate}
\end{theorem}

\begin{theorem} \label{th2}
\begin{enumerate}
\item We have a commutative semigroup for  $(\mathbb{R}_\digamma, \oplus)$ with the zero element $(\mathbb{R}_\digamma, \oplus)$.
\item There is no opposite element if there are fuzzy numbers which are not crisp ($(\mathbb{R}_\digamma, \oplus)$ cannot be a group).
\item $ \forall z_1, z_2 \in \mathbb{R}$ with $z_1, z_2 \geq 0$ or $z_1, z_2 \leq 0$ and  $ \forall p \in \mathbb{R}_\digamma$, one get $(z_1+z_2) \odot p = z_1 \odot p \oplus z_2 \odot u$.
%For arbitrary $z_1, z_2 \in \mathbb{R}$, this property is not fulfilled.
\item  $\forall \gamma \in \mathbb{R}$ and $p, p_1 \in \mathbb{R}_\digamma$, one get $\gamma \odot (p \oplus p_1) = \gamma \odot p \oplus \gamma \odot p_1 $
\item   $ \forall \gamma, \mu \in \mathbb{R}$ and $p \in \mathbb{R}_\digamma$, one get $\gamma \odot (\mu \odot p) = (\gamma \mu) \odot p$.
\item There is the general attributes of the norm for of $\| . \|_\digamma: \mathbb{R}_\digamma \rightarrow \mathbb{R}$ by $\| p\|_\digamma = \mathcal{D}(p,\tilde{0})$ which is $\| p\|_\digamma =0  \Leftrightarrow  p  = \tilde{0}$,
    $
    \| \gamma \odot p \|_\digamma =  | \gamma | \| p\|_\digamma
    $ and $\| p  \oplus p_1 \|_\digamma \leq  \| p\|_\digamma + \| p_1\|_\digamma$
\item $| \| p\|_\digamma + \| p_1\|_\digamma | \leq \mathcal{D}(p,p_1)$ and $\mathcal{D}(p,p_1) \leq \| p\|_\digamma + \| p_1\|_\digamma$ for any $p,p_1 \in \mathbb{R}_\digamma$.
\end{enumerate}
\end{theorem}

\begin{definition} \label{def5} Continuity of a fuzzy real number valued function $Y: [z_1,z_2] \rightarrow \mathbb{R}_\digamma$ can be defined in $x_0 \in [z_1,z_2]$ as
$\forall \varepsilon > 0, ~\exists \rho > 0;~\mathcal{D}(Y(x), Y(x_0)) <  \varepsilon $, whenever $x \in [z_1,z_2]$ and $| x-x_0 | < \rho$. %We say that $Y$ is fuzzy continuous on $[z_1,z_2]$ if $Y$ is continuous at each $x_0 \in [z_1,z_2]$ and denote the space of all such functions by $C_\digamma[z_1,z_2]$.
\end{definition}

\begin{definition} \label{def6}
Assume that  $Y: [z_1,z_2] \rightarrow \mathbb{R}_\digamma$ is a bounded mapping. The modulus of continuity $\omega_{[z_1,z_2]}(Y,.): \mathbb{R}_+ \cup \{ 0 \}  \rightarrow \mathbb{R}_+$ is defined as
\begin{equation}\label{1}
\omega_{[z_1,z_2]}(Y,\rho) =  \sup \{ \mathcal{D}(Y(x), Y(y)): x,y \in [z_1,z_2], |x-y| \leq \rho \}.
\end{equation}
%where $\mathbb{R}_+$ shows the set of positive real numbers, is called the modulus of oscillation of Y on $[z_1,z_2]$.
Also $\omega_{[z_1,z_2]}(Y,\rho)$ is the uniform modulus of continuity of $Y$ if $Y \in C_\digamma [z_1,z_2]$.
\end{definition}
%Some properties of the modulus of continuity are given below:

\begin{theorem} \label{th2}
We have the following properties for the modulus of continuity:
\begin{enumerate}
  \item $\mathcal{D}(Y(x), Y(y))  \leq \omega_{[z_1,z_2]}(Y, | x-y |)$ for any $x,y \in [z_1,z_2]$.
  \item $\omega_{[z_1,z_2]}(Y,\rho)$ is increasing function of $\rho$,
  \item $\omega_{[z_1,z_2]}(Y,0) = 0$,
  \item  $\omega_{[z_1,z_2]}(Y,\rho_1+\rho_2)  \leq  \omega_{[z_1,z_2]}(Y,\rho_1) + \omega_{[z_1,z_2]}(Y,\rho_2), ~\rho_1, \rho_2 \geq 0 $
  \item $\omega_{[z_1,z_2]}(Y, n \rho) \leq n \omega_{[z_1,z_2]}(Y,\rho)$  for any $\rho \geq 0$ and $n \in N$,
  \item  $\omega_{[z_1,z_2]}(Y, \gamma \rho) \leq (\gamma +1) \omega_{[z_1,z_2]}(Y, \rho),~ \rho, \gamma \geq 0$,
  \item For $[z_3,z_4] \subseteq [z_1,z_2]$ one get  $\omega_{[z_3,z_4]}(Y,\rho) \leq \omega_{[z_1,z_2]}(Y,\rho)$.
\end{enumerate}
\end{theorem}

\begin{definition} \label{def7}
Assume that $Y: [z_1,z_2] \rightarrow \mathbb{R}_\digamma$. $Y $ is a Riemann integrable of fuzzy type to $I(Y) \in \mathbb{R}_\digamma$ if  $ \forall \varepsilon > 0,~ \exists \rho > 0$; $ \forall$ division $P = \{ [p,p_1]: \xi \}$ of $[z_1,z_2]$ with the norms $\Delta (p) < \rho$, it holds
\begin{equation}\label{2}
\mathcal{D} \left( {\sum_p}^* (p_1-p) \odot Y (\xi), I(Y)   \right) < \varepsilon;
\end{equation}
 where  $\sum^*$ shows the fuzzy summation. Then 
$$
I(Y) = (\mathcal{FR}) \int_{z_1}^{z_2} Y(x) dx.
$$
And for $Y \in C_\digamma[z_1,z_2]$ it follows
$$
\begin{array}{l}
  \underline{(\mathcal{FR}) \int_{z_1}^{z_2} Y(t;r) dt} = \int_{z_1}^{z_2}\underline{Y}(t;r) dt,  \\
  \\
  \overline{(\mathcal{FR}) \int_{z_1}^{z_2} Y(t;r) dt} = \int_{z_1}^{z_2}\overline{Y}(t;r) dt
\end{array}
$$
\end{definition}

\begin{lemma} \label{l1}
If $Y,V: [z_1,z_2] \subseteq \mathbb{R} \rightarrow \mathbb{R}_\digamma$ are fuzzy and continuous functions, then $Y: [z_1,z_2] \rightarrow \mathbb{R}_+$ by $F(x) = \mathcal{D}(Y(x), V(x))$ is
continuous on $[z_1,z_2]$ and
\begin{equation}\label{3}
\mathcal{D} \left( (\mathcal{FR}) \int_{z_1}^{z_2} Y(x) dx, (\mathcal{FR}) \int_{z_1}^{z_2} V(x) dx   \right) \leq \int_{z_1}^{z_2} \mathcal{D}(Y(x), V(x)) dx.
\end{equation}
\end{lemma}

\begin{theorem}\label{th4}
  Assume that  $Y: [z_1,z_2] \rightarrow \mathbb{R}_\digamma$ is a Henstock integrable and a bounded function. Then for
 $z_1= x_0 < x_1 < ... < x_n = z_2$ and $\xi_i \in [x_{i-1},x_i]$ it gives:
 $$
 \mathcal{D} \left(  (\mathcal{FH}) \int_{z_1}^{z_2} Y(t) dt, \sum_{i=1}^{n}~^* (x_i - x_{i-1}) \odot Y(\xi_i)  \right)  \leq \sum_{i=1}^n (x_i - x_{i-1}) \omega_{[x_i, x_{i-1}]} (Y, x_i - x_{i-1})
  $$
\end{theorem}

\begin{corollary}\label{co1}
 Let $Y: [z_1,z_2] \rightarrow \mathbb{R}_\digamma$ be a bounded and Henstock integrable function. Then
\begin{enumerate}
  \item $ \mathcal{D} \left( (\mathcal{FH}) \int_{z_1}^{z_2} Y(t) dt,  (z_2-z_1) \odot Y(\frac{z_1+z_2}{2}) \right) \leq \frac{z_2-z_1}{2} \omega_{[z_1,z_2]} (Y, \frac{z_2-z_1}{2} ) $
  \item $\mathcal{D} \left( (\mathcal{FH}) \int_{z_1}^{z_2} Y(t) dt,  \frac{z_2-z_1}{2} \odot (Y(z_1) \oplus Y(z_2)) \right) \leq \frac{z_2-z_1}{2} \omega_{[z_1,z_2]} (Y, \frac{z_2-z_1}{2} ) $
  \item $\mathcal{D} \left( (\mathcal{FH}) \int_{z_1}^{z_2} Y(t) dt,  \frac{z_2-z_1}{6} \odot (Y(z_1) \oplus 4 \odot Y(\frac{z_1+z_2}{2}) \oplus Y(z_2)) \right) \leq  2 (z_2-z_1) \omega_{[z_1,z_2]} (Y, \frac{z_2-z_1}{6} ). $
\end{enumerate}
 \end{corollary}

%Since any Riemann integrable fuzzy number valued function is also Henstock integrable, it follows that the
%above quadrature rules hold for Riemann integrable too.

%%%======================================================================================================================
%%%======================================================================================================================
%%%======================================================================================================================
\section{Main Idea}

In this section let us discuss the  existence and uniqueness of the solution of problem (\ref{4}) based on the successive approximations.  Assume  $X = \{ Y:[z_1,z_2] \rightarrow \mathbb{R}_\digamma: Y ~is ~ continuous \}$ is the continuous functions space with fuzzy distance $\mathcal{D}^* (Y,V) = \sup_{z_1 \leq v  \leq z_2} \mathcal{D}(Y(v), V(v))$. Let $A: X \rightarrow X$ be a nonlinear integral operator. Application of $A$ for the problem (\ref{4}) gives
$$
A Z(v) = Y(v) \oplus (\mathcal{FR}) \sum_{v=1}^{m'} \int_{\theta_{t-1}(v)}^{\theta_{t}(v)} K_t(r,v) \odot G (Z(r)) dr,~~~ \forall s,v \in [z_1,z_2],  \forall F \in X.
$$
Then 

\begin{theorem} \label{th5}
Assume that the kernels  $K_1(r,v), K_2(r,v), ..., K_{m'}(r,v),~z_1 \leq s,v \leq T \leq z_2$ are positive and continuous. Let function $Y(v)$ be a  fuzzy continuous of $v$, $z_1 \leq v \leq T \leq z_2$. Moreover 
$$
\exists L>0; ~ \mathcal{D}(G(Z_1(p)), G(Z_2(p_1))) \leq L \mathcal{D} (Z_1(p), Z_2(p_1)),~~~ \forall p,p_1 \in [z_1,z_2].
$$
If $c = \sum_{t=1}^{m'} M_t L (\theta_{t} - \theta_{t-1}) < 1$ then there is a unique solution $F^* \in X$ for  the FVIE (\ref{4}) based on the following successive approximations method:
 \begin{equation}\label{5}\left\{
\begin{array}{l}
\displaystyle  Z_0(v)=Y(v), \\
   \\
\displaystyle Z_m(v) = Y(v) \oplus (\mathcal{FR}) \sum_{t=1}^{m'} \int_{\theta_{t-1}(v)}^{\theta_{t}(v)} K_t(r,v) \odot G (Z_{m-1}(r)) dr,~~~z_1 \leq r,v \leq T \leq z_2, ~~m \geq 1,
\end{array}\right.
\end{equation}
which is convergent to $F^*$. Also 
\begin{equation}\label{6}
\mathcal{D}(F^*(v), Z_m(v)) \leq\frac{c^{m+1}}{L (1-c)} M_0, ~~~\forall t \in [z_1,z_2],~~m\geq  1
\end{equation}
is the error bound for $M_0 = \sup_{z_1 \leq v  \leq z_2}   \| G(Y(v))\|_\digamma $.
\end{theorem}

\textbf{Proof: } We use the Banach fixed point principle to prove the theorem. Let us show $A: X \rightarrow X$ and also prove the uniformly continuity of the operator $A$. We know the continuity of $Z$ on the compact set of $[z_1,z_2]$ thus that is uniformly continuous and it follows
$$
\forall \varepsilon_1 > 0 \exists\rho_1 >0; \mathcal{D}(Z(v_1), Z(v_2)) < \varepsilon_1~ whenever~ | v_1-v_2 | < \rho_1, \forall v_1, v_2 \in [z_1,z_2].
$$

Also  $K_t, t=1,2, ..., m'$ is uniformly continuous. Therefore for $\varepsilon_t >0$ there is an estimate
$$
| K_t(r,v_1) - K_t(r,v_2) | < \varepsilon_t~ whenever ~  | v_1 - v_2| < \rho_t,~ \forall v_1, v_2 \in [z_1,z_2].
$$

Assume that $\rho = \min \{\rho_1, \rho_2, ..., \rho_{m'}\}$ and  $v_1, v_2 \in [z_1,z_2]$ with $| v_1 - v_2| < \rho_t$.  Applying Lemma \ref{l1}  and Theorem \ref{th1} one can write:

$$
\begin{array}{l}
  \mathcal{D} (A(F)(v_1),A(F)(v_2) ) \\
   \\
 \displaystyle  \leq  \mathcal{D} (Y(v_1), Y(v_2)) + \mathcal{D} ((\mathcal{FR}) \sum_{t=1}^{m'}\int_{\theta_{t-1}(v_1)}^{\theta_{t}(v_1)} K_t(r,v_1) \odot G (Z(r)) dr, (\mathcal{FR}) \sum_{t=1}^{m'}\int_{\theta_{t-1}(v_2)}^{\theta_{t}(v_2)} K_t(r,v_2) \odot G (Z(r)) dr)\\
    \\
 \displaystyle  \leq \varepsilon_1 + L  \sum_{t=1}^{m'}  | K_t(r,v_1) - K_t(r,v_2) | \mathcal{D} ((\mathcal{FR})  \int_{\theta_{t-1}(v_1)}^{\theta_{t}(v_1)} G (Z(r)) dr, (\mathcal{FR})  \int_{\theta_{t-1}(v_2)}^{\theta_{t}(v_2)} \tilde{0} dr)\\
  \\
  \displaystyle \leq \varepsilon_1  + \sum_{t=1}^{m'} \varepsilon_t (\theta_{t}(v_1) - \theta_{t-1}(v_1)) M_0
\end{array}
$$
where $M_0 = \sup_{z_1 \leq s \leq T \leq z_2} \| G(Y(r))  \|_\digamma$. By choosing $\varepsilon_1 = \frac{\varepsilon}{m'+1}$ and $\varepsilon_t = \frac{\varepsilon}{(m'+1) (\theta_{t}(v_1) - \theta_{t-1}(v_1))M_0}$ we find $\mathcal{D} (A(F)(v_1),A(F)(v_2) ) < \varepsilon$.

Thus $A(F)$ is uniformly continuous for any $F \in X$, and so continuous on $[z_1,z_2]$, and hence $A(X) \subset X$. Now, it can be proved  the contracting map of
 the operator $A$. For $Z_1, Z_2 \in X$ and $t \in [z_1,z_2]$
$$
\begin{array}{l}
  \mathcal{D} (A(Z_1)(v),A(Z_2)(v) ) \\
   \\
 \displaystyle  \leq  \mathcal{D} (Y(v), Y(v)) + \mathcal{D} ((\mathcal{FR}) \sum_{t=1}^{m'}\int_{\theta_{t-1}(v)}^{\theta_{t}(v)} K_t(r,v) \odot G (Z_1(r)) dr, (\mathcal{FR}) \sum_{t=1}^{m'}\int_{\theta_{t-1}(v)}^{\theta_{t}(v)} K_t(r,v) \odot G (Z_2(r)) dr)\\
    \\
 \displaystyle  \leq   \sum_{t=1}^{m'} \int_{\theta_{t-1}(v)}^{\theta_{t}(v)}  \mathcal{D} (  K_t(r,v) \odot G (Z_1(r)) , (\mathcal{FR}) K_t(r,v) \odot G (Z_2(r))) dr\\
  \\
  \displaystyle \leq L \sum_{t=1}^{m'} M_t (\theta_{t}(v) - \theta_{t-1}(v)) \mathcal{D}^*(Z_1,Z_2) = C \mathcal{D}^* (Z_1,Z_2),
\end{array}
$$
thus, $\mathcal{D} (A(Z_1)(v),A(Z_2)(v) )\leq C \mathcal{D}^* (Z_1,Z_2)$. As  $C < 1$ and $A$ is a contraction on the Banach space $(X,\mathcal{D}^*)$. Thus based on the Banach fixed point principle there is unique solution $F^*$ in $X$ for Eq. (\ref{4}) and
$$
\mathcal{D}(F^*(v), Z_m(v)) \leq \mathcal{D}^*(F^*,Z_m) \leq  \frac{C^m}{1-C} \mathcal{D}^*(Z_0,Z_1), ~~z_1 \leq v \leq T \leq z_2,~~m\geq  1.
$$
Also one can write
$$
\begin{array}{l}
\displaystyle  \mathcal{D}^*(Z_0,Z_1) =  \sup_{z_1 \leq v  \leq z_2} \mathcal{D}(Y(v), Y(v)  +  (\mathcal{FR}) \sum_{t=1}^{m'}\int_{\theta_{t-1}(v)}^{\theta_{t}(v)} K_t(r,v) \odot G (Z(r)) dr )\\
    \\
\displaystyle  \leq  \sup_{z_1 \leq v  \leq z_2}  \sum_{t=1}^{m'}\int_{\theta_{t-1}(v)}^{\theta_{t}(v)} \mathcal{D}(\tilde{0},K_t(r,v) \odot G (Z_0(r)) ) dr\\
  \\
\displaystyle  \leq    \sum_{t=1}^{m'} M_t \int_{\theta_{t-1}(v)}^{\theta_{t}(v)}  \sup_{z_1 \leq v  \leq z_2} \mathcal{D}(\tilde{0}, G (Z_0(r)) ) dr \\
  \\
\displaystyle  = \sum_{t=1}^{m'} M_t (\theta_{t}(v) - \theta_{t-1}(v)) M_0 = \frac{C}{L} M_0.
\end{array}
$$

Now it can be introducde the following numerical method to find the approximate solution of (\ref{4}). As
$$
z_1=v_0 < v_1 < ... < v_{n-1} < v_n = z_2
$$
where $v_i = a+ih$ and $h=\frac{b-a}{n}$ and one have the following iterative procedure as
 $$\left\{
\begin{array}{l}
\displaystyle  y_0(v)=Y(v), \\
   \\
\displaystyle y_m(v) = Y(v) \oplus  \sum_{t=1}^{m'} \frac{h}{2}\odot \bigg[ K_t(v_0,v) \odot G(y_{m-1}(v_0)) \oplus K_t(v_n,v) \odot G(y_{m-1}(v_n)) \\
\\
~~~~~~~~~~~~~~~~~~~~~~~~~~~~~~~~~~~\displaystyle \oplus 2 \sum_{l=1}^{n-1} K_t(v_l,v) \odot G(y_{m-1}(v_l)) \bigg], ~~m \geq 1.
\end{array}\right.
$$
Also the compact form of the relation is
\begin{equation}\label{8}
 \left\{
\begin{array}{l}
\displaystyle  y_0(v)=Y(v), \\
   \\
\displaystyle y_m(v) = Y(v) \oplus  \sum_{t=1}^{m'}  \sum_{l=1}^{n-1} \frac{h}{2}\odot \bigg[ K_t(v_l,v) \odot G(y_{m-1}(v_l)) \oplus K_t(v_l,v) \odot G(y_{m-1}(v_l))  \bigg], ~~m \geq 1.
\end{array}\right.
\end{equation}

%%%======================================================================================================================
%%%======================================================================================================================
%%%======================================================================================================================
\section{Error Estimation}

\begin{theorem} \label{th6}
Assume that the nonlinear FVIE (\ref{4}) with kernel $K_t(r,v)$ along continuous curves $\theta_t(v), t=1,2,...,m'$ with positive sign on $[z_1,z_2] \times [z_1,z_2]$, $G$ continuous on $\mathbb{R}_\digamma$ and $Y$ continuous on $[z_1,z_2]$. Moreover there exists $L > 0$ such that
$$
\mathcal{D}(G(Z_1(p)),G(Z_2(p_1)) ) \leq L. \mathcal{D} (Z_1(p),Z_2(p_1) ), ~~\forall p,p_1 \in [z_1,z_2].
$$
For $C_t = M_t L (z_2-z_1) < 1 $ where $M_t = \max_{z_1 \leq r,v \leq T \leq z_2}  | K_t(r,v)|$, then the successive scheme (\ref{8}) converges to the unique solution of (\ref{4}), $F$ and the error estimation can be obtained as:
$$
\mathcal{D}^* (F,y_m) \leq \sum_{t=1}^{m'} \frac{C_t}{2(1-C_t)} \omega_{\theta_{t-1},\theta_t}(Y,h) + \sum_{t=1}^{m'} \frac{C_t^{m+1} L_1 }{L(1-C_t)} + \sum_{t=1}^{m'} \frac{C_t^2+2C_t}{2 L M_t(1-C_t)} (L_1\omega_s(K_t,h)+L_2\omega_t(K_t,h) )
$$
where
$$\omega_s(K_t,h) = \sup_{z_1 \leq v \leq T \leq z_2} \{  \sup | K_t(x,v) - K_t(y,v) |: | x-y  | \leq h  \},~t=1,2,...,m', $$
and
$$\omega_t(K_t,h) = \sup_{z_1 \leq s \leq T \leq z_2} \{  \sup | K_t(r,v_1) - K_t(r,v_2) |: | v_1-v_2  | \leq h  \},~t=1,2,...,m'. $$
\end{theorem}

\textbf{Proof:}
We know
$$
Z_1(v) = Y(v) \oplus (\mathcal{FR}) \sum_{t=1}^{m'} \int_{\theta_{t-1}(v)}^{\theta_{t}(v)} K_t(r,v) \odot G (Z_0(r)) dr,~~~z_1 \leq s,v \leq T \leq z_2,
$$
then
$$
\begin{array}{l}
\displaystyle \mathcal{D}(Z_1(v), y_1(v)) = \mathcal{D}(Y(v),Y(v)) \\
\\
\displaystyle + \mathcal{D} \bigg((\mathcal{FR})\sum_{t=1}^{m'} \int_{\theta_{t-1}(v)}^{\theta_{t}(v)} K_t(r,v) \odot G (Z_0(r)) dr, \\
\\
\displaystyle ~~~~~~~~ \sum_{t=1}^{m'} \sum_{l=0}^{n-1} \frac{h}{2}\odot \big[ K_t(v_l,v)\odot G(Z_0(v_t)) \oplus K_t (v_l,v) \odot G(Z_0(v_l)) \big] \bigg)    \\
    \\
\displaystyle = \mathcal{D} \bigg(\sum_{l=0}^{n-1} \sum_{t=1}^{m'} (\mathcal{FR}) \int_{\theta_{t-1}(v_l)}^{\theta_{t}(v_l)} K_t(r,v) \odot G (Y(r)) dr,   \\
\\
\displaystyle ~~~~~~~~ \sum_{t=1}^{m'} \sum_{l=0}^{n-1} \frac{h}{2}\odot \big[ K_t(v_l,v)\odot G(Y(v_l)) \oplus K_t (v_{l+1},v) \odot G(Y(v_{l+1})) \big] \bigg)\\
\\
\leq
\displaystyle \sum_{t=1}^{m'} \sum_{l=0}^{n-1} \mathcal{D} \bigg(  (\mathcal{FR}) \int_{\theta_{t-1}(v_l)}^{\theta_{t}(v_{l+1})} K_t(r,v) \odot G (Y(r)) dr, \\
\\
\displaystyle  ~~~~~~~~ \frac{h}{2}\odot \big[ K_t(v_l,v)\odot G(Y(v_l)) \oplus K_t (v_{l+1},v) \odot G(Y(v_{l+1})) \big] \bigg)\\
\\
\leq
\displaystyle \sum_{t=1}^{m'} \sum_{l=0}^{n-1} \mathcal{D} \bigg(  (\mathcal{FR}) \int_{\theta_{t-1}(v_l)}^{\theta_{t}(v_{l+1})} K_t(r,v) \odot G (Y(r)) dr,  \\
\\
\displaystyle ~~~~~~~~ \frac{h}{2}\odot \big[ K_t(r,v)\odot G(Y(v_l)) \oplus K_t (r,v) \odot G(Y(v_{l+1})) \big] \bigg) \\
  \\
\displaystyle  ~~~~~~~~+  \sum_{t=1}^{m'} \sum_{l=0}^{n-1} \mathcal{D} \bigg(    \frac{h}{2}\odot \big[ K_t(r,v)\odot G(Y(v_l)) \oplus K_t (r,v) \odot G(Y(v_{l+1})) \big], \\
  \\
\displaystyle ~~~~~~~~ \frac{h}{2}\odot \big[ K_t(v_l,v)\odot G(Y(v_l)) \oplus K_t (v_{l+1},v) \odot G(Y(v_{l+1})) \big]   \bigg)
\end{array}
$$

Applying the second part of the first corollary and regarding to Lemma 4 in [11] one have:
$$
\begin{array}{ll}
\displaystyle \mathcal{D}(Z_1(v), y_1(v)) &\leq
\displaystyle \sum_{t=1}^{m'} \sum_{l=0}^{n-1}  |K_t(r,v)|   \mathcal{D} \bigg(  (\mathcal{FR}) \int_{\theta_{t-1}(v_l)}^{\theta_{t}(v_{l+1})}   G (Y(r)) dr,  \frac{h}{2}\odot \big[  G(Y(v_l)) \oplus  G(Y(v_{l+1})) \big] \bigg) \\
  \\
&\displaystyle  ~~~+  \sum_{t=1}^{m'} \sum_{l=0}^{n-1} \mathcal{D} \bigg(    \frac{h}{2}\odot  K_t(r,v)\odot G(Y(v_l))  ,  \frac{h}{2}\odot K_t(v_l,v)\odot G(Y(v_l))   \bigg)\\
\\
&\displaystyle  ~~~+  \sum_{t=1}^{m'} \sum_{l=0}^{n-1} \mathcal{D} \bigg(    \frac{h}{2}\odot  K_t(r,v)\odot G(Y(v_l))  , \frac{h}{2}\odot K_t(v_{l+1},v)\odot G(Y(v_{l+1}))   \bigg)\\
\\
&\displaystyle \leq \frac{h}{2} \sum_{t=1}^{m'} \sum_{l=0}^{n-1} |K_t(r,v)|  \omega_{[v_l,v_{l+1}]} (G(Y),\frac{h}{2} ) \\
\\
&\displaystyle~~~+ \frac{h}{2} \sum_{t=1}^{m'} \sum_{l=0}^{n-1}   |K_t(r,v) -K_t(v_{r},v) | \mathcal{D}(G(Y(v_l)),\tilde{0}) \\
\\
&\displaystyle ~~~ +  \frac{h}{2} \sum_{t=1}^{m'} \sum_{l=0}^{n-1}   |K_t(r,v) -K_t(v_{l+1},v) | \mathcal{D}(G(Y(v_{l+1})),\tilde{0})\\
\\
&\displaystyle \leq \sum_{t=1}^{m'} \frac{M_t(\theta_t - \theta_{t-1})}{2} \omega_{[ \theta_{t-1},\theta_t]} (G(Y),h) + \sum_{t=1}^{m'} (\theta_t - \theta_{t-1}) M_0\omega_s (K_t,h)\\
\\
&\displaystyle = \sum_{t=1}^{m'} \frac{M_t(\theta_t - \theta_{t-1})}{2}  \sup_{\theta_{t-1} \leq p,p_1 \leq \theta_t}   \{ \mathcal{D}(G(Y(p)), G(Y(p_1))): |u-p_1  | \leq h  \} \\
\\ & \displaystyle ~~~ +  \sum_{t=1}^{m'}  (\theta_t - \theta_{t-1}) M_0  \omega_s (K_t,h)\\
\\
&\displaystyle \leq \sum_{t=1}^{m'} \frac{M_t(\theta_t - \theta_{t-1})}{2}  \sup_{\theta_{t-1} \leq p,p_1 \leq \theta_t}   \{ L . \mathcal{D}(Y(p),Y(p_1)): |u-p_1  | \leq h  \} \\
\\ & \displaystyle ~~~ +  \sum_{t=1}^{m'}  (\theta_t - \theta_{t-1}) M_0  \omega_s (K_t,h)\\
\\
&\displaystyle \leq \sum_{t=1}^{m'} \frac{M_t L (\theta_t - \theta_{t-1})}{2}  \omega_{[\theta_{t-1}, \theta_t]}   (Y,h) +  \sum_{t=1}^{m'}  (\theta_t - \theta_{t-1}) M_0  \omega_s (K_t,h)\\
\end{array}
$$
where $M_0 = \sup_{\theta_{t-1} \leq s \leq \theta_t}  \| G(Y(r))\|_\digamma$ and $\omega_s (K_t,h)$ is the partial modules of continuity. Thus
$$
\mathcal{D}(Z_1(v), y_1(v)) \leq \sum_{t=1}^{m'} \frac{C_t}{2}  \omega_{[\theta_{t-1}, \theta_t]}   (Y,h) +  \sum_{t=1}^{m'}  \frac{C_t}{L M_t} M_0  \omega_s (K_t,h).
$$
We have
$$
Z_2(v) = Y(v) \oplus (\mathcal{FR}) \sum_{t=1}^{m'} \int_{\theta_{t-1}(v)}^{\theta_{t}(v)} K_t(r,v) \odot G (Z_1(r)) dr,
$$
therefore
$$
\begin{array}{l}
\displaystyle \mathcal{D}(Z_2(v), y_2(v)) =
\displaystyle \mathcal{D}(Y(v),Y(v)) \\
\\
+   \mathcal{D} \bigg(  (\mathcal{FR}) \sum_{t=1}^{m'} \int_{\theta_{t-1}(v)}^{\theta_{t}(v)} K_t(r,v) \odot  G (Z_1(r)) dr, \\
\\
\displaystyle ~~~~~~~\sum_{t=1}^{m'} \sum_{l=0}^{n-1} \frac{h}{2}\odot \big[ K_t(v_l,v) \odot G(y_1(v_l)) \oplus K_t(v_{l+1},v) \odot G(y_1(v_{l+1}))  \big] \bigg) \\
\\
 \displaystyle \leq \mathcal{D} \bigg(  (\mathcal{FR}) \sum_{t=1}^{m'} \int_{\theta_{t-1}(v)}^{\theta_{t}(v)} K_t(r,v) \odot  G (Z_1(r)) dr, \\
 \\
\displaystyle ~~~~~~~ \sum_{t=1}^{m'} \sum_{l=0}^{n-1} \frac{h}{2}\odot \big[ K_t(r,v) \odot G(Z_1(v_l)) \oplus K_t(r,v) \odot G(Z_1(v_{l+1}))  \big] \bigg)\\
\\
 \displaystyle + \mathcal{D} \bigg(  \sum_{t=1}^{m'} \sum_{l=0}^{n-1} \frac{h}{2}\odot \big[ K_t(r,v) \odot G(Z_1(v_l)) \oplus K_t(r,v) \odot G(Z_1(v_{l+1}))  \big], \\
 \\
 \displaystyle ~~~~~~~ \sum_{t=1}^{m'} \sum_{l=0}^{n-1} \frac{h}{2}\odot \big[ K_t(r,v) \odot G(y_1(v_l)) \oplus K_t(r,v) \odot G(y_1(v_{l+1}))  \big] \bigg)\\
 \\
  \displaystyle + \mathcal{D} \bigg(  \sum_{t=1}^{m'} \sum_{l=0}^{n-1} \frac{h}{2}\odot \big[ K_t(r,v) \odot G(y_1(v_l)) \oplus K_t(r,v) \odot G(y_1(v_{l+1}))  \big], \\
 \\
 \displaystyle ~~~~~~~ \sum_{t=1}^{m'} \sum_{l=0}^{n-1} \frac{h}{2}\odot \big[ K_t(v_l,v) \odot G(y_1(v_l)) \oplus K_t(v_{l+1},v) \odot G(y_1(v_{l+1}))  \big] \bigg)\\
 \\
\displaystyle \leq \sum_{t=1}^{m'} \frac{M_t(\theta_t - \theta_{t-1})}{2} \omega_{[\theta_{t-1}, \theta_{t}]} (G(Z_1), \frac{h}{2})  + \sum_{t=1}^{m'} \frac{M_t(\theta_t - \theta_{t-1})}{2} \bigg[ \mathcal{D}(G(Z_1(v_l)), G(y_1(v_l)))\\
\\
\displaystyle ~~~~~~~ + \mathcal{D}(G(Z_1(v_{l+1})), G(y_1(v_{l+1}))) \bigg] + \sum_{t=1}^{m'}(\theta_t - \theta_{t-1}) M_1 \omega_s(K_t,h)\\
\\
\displaystyle \leq \sum_{t=1}^{m'} \frac{C_t}{2} \omega_{[\theta_{t-1}, \theta_{t}]} (Z_1, h)  + \sum_{t=1}^{m'} \frac{C_t}{2} \bigg[ \mathcal{D}(Z_1(v_l), y_1(v_l)) + \mathcal{D}(Z_1(v_{l+1}), y_1(v_{l+1})) \bigg] + \sum_{t=1}^{m'} \frac{C_t}{L M_t} M_1 \omega_s(K_t,h)
\end{array}
$$
where $M_1 = \sup_{\theta_{t-1} \leq s \leq \theta_t}  \| G(y_1(r))\|_\digamma$. Applying (\ref{5}) and (\ref{8}) and using induction for $m\geq 3$ it can be concluded
\begin{equation}\label{9}
\begin{array}{ll}
\displaystyle \mathcal{D}(Z_m(v), y_m(v)) & \displaystyle \leq \sum_{t=1}^{m'} \frac{C_t}{2} \omega_{[\theta_{t-1}, \theta_{t}]} (Z_{m-1}, h) + \sum_{t=1}^{m'} \frac{C_t}{2} \bigg[ \mathcal{D}(Z_{m-1}(v_l), y_{m-1}(v_l)) \\
\\
& \displaystyle  + \mathcal{D}(Z_{m-1}(v_{l+1}), y_{m-1}(v_{l+1})) \bigg] + \sum_{t=1}^{m'} \frac{C_t}{L M_t} M_{m-1} \omega_s(K_t,h)
\end{array}
\end{equation}
where $M_{m-1} = \sup_{\theta_{t-1} \leq v \leq \theta_t}  \| G(y_{m-1}(v))\|_\digamma$. Taking supremum for $z_1 \leq v \leq z_2 $ from (\ref{9}) we have
\begin{equation}\label{10}
\begin{array}{l}
\displaystyle \mathcal{D}^*(Z_m, y_m) \leq \sum_{t=1}^{m'} \frac{C_t}{2} \omega_{[\theta_{t-1}, \theta_{t}]} (Z_{m-1}, h)  + \sum_{t=1}^{m'} C_t \displaystyle\mathcal{D}^*(Z_{m-1}, y_{m-1}) + \sum_{t=1}^{m'} \frac{C_t}{L M_t} M_{m-1} \omega_s(K_t,h),\\
\\
\displaystyle \mathcal{D}^*(Z_{m-1}, y_{m-1}) \leq \sum_{t=1}^{m'} \frac{C_t}{2} \omega_{[\theta_{t-1}, \theta_{t}]} (Z_{m-2}, h)  + \sum_{t=1}^{m'} C_t \displaystyle\mathcal{D}^*(Z_{m-2}, y_{m-2}) + \sum_{t=1}^{m'} \frac{C_t}{L M_t} M_{m-2} \omega_s(K_t,h),\\
\\
\vdots\\
\\
\displaystyle \mathcal{D}^*(Z_1, y_1) \leq \sum_{t=1}^{m'} \frac{C_t}{2} \omega_{[\theta_{t-1}, \theta_{t}]} (Y, h)   + \sum_{t=1}^{m'} \frac{C_t}{L M_t} M_{0} \omega_s(K_t,h).
\end{array}
\end{equation}
If one multiple the above inequality to $1, C_t, ..., C_t^{m-1}$ and find the summation then
\begin{equation}\label{11}
\begin{array}{ll}
\displaystyle \mathcal{D}^*(Z_m, y_m) & \displaystyle \leq \sum_{t=1}^{m'} \frac{C_t}{2} \bigg( \omega_{[\theta_{t-1}, \theta_{t}]} (Z_{m-1}, h)  +  C_t \omega_{[\theta_{t-1}, \theta_{t}]} (Z_{m-2}, h) + ... + C_t^{m-1} \omega_{[\theta_{t-1}, \theta_{t}]} (Y, h) \bigg)\\
\\
&\displaystyle + \sum_{t=1}^{m'} \frac{C_t}{L M_t} \omega_s(K_t,h) (M_{m-1}+ C_t M_{m-2} + ... + C_t^{m-1} M_0).\\
\end{array}
\end{equation}
Moreover for $v_1, v_2 \in [z_1,z_2]$ with $| v_1 - v_2 \leq h  |$  one can write
$$
\begin{array}{l}
  \mathcal{D} (Z_m(v_1),Z_m(v_2) ) \\
   \\
 \displaystyle  =   \mathcal{D} \bigg( Y(v_1) \oplus \odot (\mathcal{FR}) \sum_{t=1}^{m'}\int_{\theta_{t-1}(v_1)}^{\theta_{t}(v_1)} K_t(r,v_1) \odot G (Z_{m-1}(r)) dr, \\
 \\
 \displaystyle ~~~~~~~~~ Y(v_2) \oplus \odot (\mathcal{FR}) \sum_{t=1}^{m'}\int_{\theta_{t-1}(v_2)}^{\theta_{t}(v_2)} K_t(r,v_2) \odot G (Z_{m-1}(r)) dr \bigg)\\
    \\
 \displaystyle  \leq  \mathcal{D} ( Y(v_1), Y(v_2)) +     \sum_{t=1}^{m'}   \int_{\theta_{t-1}(v_1)}^{\theta_{t}(v_1)} | K_t(r,v_1) - K_t(r,v_2) |  \mathcal{D}( G (Z_{m-1}(r)), \tilde{0} ) dr  \\
 \\
 \displaystyle \leq  \mathcal{D} ( Y(v_1), Y(v_2)) +     \sum_{t=1}^{m'}  \frac{C_t}{L M_t} \omega_t(K_t,h) M'_{m-1},
\end{array}
$$
where $\omega_t(K_t,h)$ is the partial modulus of continuity with respect to $t$. Let $M'_{m-1} =\sup_{\theta_{t-1} \leq s \leq \theta_t}  \| G(Z_{m-1}(r))\|_\digamma $ then we can find the relation $Z_m$ and $Y$ as:
\begin{equation}\label{12}
\omega_{[\theta_{t-1}, \theta_t]} (Z_m, h) \leq \omega_{[\theta_{t-1}, \theta_t]} (Y, h) + \sum_{t=1}^{m'} \frac{C_t}{L M_t} \omega_t(K_t,h) M'_{m-1}.
\end{equation}
And if we substitute above inequality into (\ref{11}) we obtain
\begin{equation}\label{13}
\begin{array}{ll}
\displaystyle \mathcal{D}^*(Z_m, y_m) & \displaystyle \leq \sum_{t=1}^{m'} \frac{C_t}{2}   (1  +  C_t + C_t^2 + ... + C_t^{m-1}) \omega_{[\theta_{t-1}, \theta_{t}]} (Y,h) \\
\\
&\displaystyle + \sum_{t=1}^{m'} \frac{C_t}{2 L M_t} \omega_t(K_t,h) ( C_t M'_{m-2} + C_t^2 M'_{m-3}+ ... + C_t^{m-1} M'_0) \\
\\
&\displaystyle + \sum_{t=1}^{m'} \frac{C_t}{ L M  } \omega_s(K_t,h) ( M_{m-1} + C_t M_{m-2} + ... + C_t^{m-1} M_0).
\end{array}
\end{equation}

Let  $L_1 = \max_{0 \leq i \leq m-1}  \{ M_i\}$ and $L_2 = \max_{0 \leq i \leq m-2}  \{ M'_i\}$ thus
\begin{equation}\label{14}
\begin{array}{ll}
\displaystyle \mathcal{D}^*(Z_m, y_m) & \displaystyle \leq \sum_{t=1}^{m'} \frac{C_t}{2}   (\frac{1-C_t^m}{1-C_t}) \omega_{[\theta_{t-1}, \theta_{t}]} (Y,h) \\
\\
&\displaystyle + \sum_{t=1}^{m'} \frac{C_t}{2 L M_t} \omega_t(K_t,h) ( C_t  + C_t^2 + ... + C_t^{m-1}) L_2 \\
\\
&\displaystyle + \sum_{t=1}^{m'} \frac{C_t}{ L M  } \omega_s(K_t,h) ( 1 + C_t + ... + C_t^{m-1} ) L_1.
\end{array}
\end{equation}

From other hand $\frac{1-C_t^m}{1-C_t} \leq \frac{1}{1-C_t}, t=1,2,...,m'$ for each $m \in N$ thus
$$
\begin{array}{ll}
\displaystyle \mathcal{D}^*(Z_m, y_m) & \displaystyle \leq \sum_{t=1}^{m'}  (\frac{ C_t}{2(1-C_t)}) \omega_{[\theta_{t-1}, \theta_{t}]} (Y,h) \\
\\
&\displaystyle + \sum_{t=1}^{m'} \frac{C_t^2+ 2 C_t}{2 L M_t (1-C_t)} (L_1 \omega_s(K_t,h) + L_2 \omega_t(K_t,h) ).
\end{array}
$$

Applying the inequality (\ref{6}) we can write
$$
\begin{array}{ll}
\displaystyle \mathcal{D}^*(F, y_m) & \leq \mathcal{D}^* (F, Z_m) + \mathcal{D}^* (Z_m, y_m)\\
\\
&\displaystyle \leq \sum_{t=1}^{m'}  (\frac{ C_t^m}{1-C_t}) \mathcal{D}^*(Z_1,Z_0) +  \sum_{t=1}^{m'}  (\frac{ C_t}{2(1-C_t)}) \omega_{[\theta_{t-1}, \theta_{t}]} (Y,h) \\
\\
&\displaystyle + \sum_{t=1}^{m'} \frac{C_t^2+ 2 C_t}{2 L M_t (1-C_t)} (L_1 \omega_s(K_t,h) + L_2 \omega_t(K_t,h) ).
\end{array}
$$
Since
$$
\begin{array}{ll}
\displaystyle \mathcal{D}(Z_1(v), Z_0(v)) &\displaystyle = \mathcal{D}\bigg( Y(v) \oplus (\mathcal{FR}) \sum_{t=1}^{m'} \int_{\theta_{t-1}(v)}^{\theta_{t}(v)} K_t(r,v) \odot G (Z_0(r)) dr, Z_0(v)\bigg) \\
\\
& \displaystyle \leq \mathcal{D}\bigg( (\mathcal{FR}) \sum_{t=1}^{m'} \int_{\theta_{t-1}(v)}^{\theta_{t}(v)} K_t(r,v) \odot G (Z_0(r)) dr, \tilde{0} \bigg)
\end{array}
$$
we obtain
$$
\begin{array}{ll}
 \displaystyle \mathcal{D}^*(Z_1, Z_0) & \displaystyle\leq  \sum_{t=1}^{m'} M_t (\theta_{t}-\theta_{t-1}) \sup_{\theta_{t-1} \leq s \leq \theta_t} \mathcal{D}(G(Y(r)), \tilde{0}) \\
    \\
  & \displaystyle = \sum_{t=1}^{m'} \frac{C_t }{L} M_0  \leq \sum_{t=1}^{m'} \frac{C_t }{L} L_1.
\end{array}
$$

Thus we get
$$
\begin{array}{ll}
\displaystyle \mathcal{D}^*(F, y_m) & \displaystyle  \leq  \sum_{t=1}^{m'}  (\frac{ C_t}{2(1-C_t)}) \omega_{[\theta_{t-1}, \theta_{t}]} (Y,h) +  \sum_{t=1}^{m'}  (\frac{ C_t^{m+1}L_1}{L(1-C_t)})    \\
\\
&\displaystyle + \sum_{t=1}^{m'} \frac{C_t^2+ 2 C_t}{2 L M_t (1-C_t)} (L_1 \omega_s(K_t,h) + L_2 \omega_t(K_t,h) ). \square
\end{array}
$$

\begin{remark}
  As we know $C_t<1, t=1,2,...,m'$ and it shows $\lim_{m  \rightarrow \infty} C_t^{m+1} = 0, t=1,2,...,m'$. And we have
  $$
  \lim_{h  \rightarrow 0} \omega_{[\theta_{t-1}, \theta_{t}]} (Y,h)=0, ~~ \lim_{h  \rightarrow 0} \omega_s (K_t,h)=0,~~  \lim_{h  \rightarrow 0} \omega_t (K_t,h)=0,~~ t=1,2,...,m'.
  $$
 The convergence of this scheme can be obtained by $\lim_{m  \rightarrow \infty, h  \rightarrow 0} \mathcal{D}^*(F,y_m)=0$. 
\end{remark}

%%%======================================================================================================================
%%%======================================================================================================================
%%%======================================================================================================================
\section{ Numerical Results}

In this section some examples are presented. We apply the mentioned method for solving the problems.

\begin{figure}[h!]
\centering
 \includegraphics[width=3in]{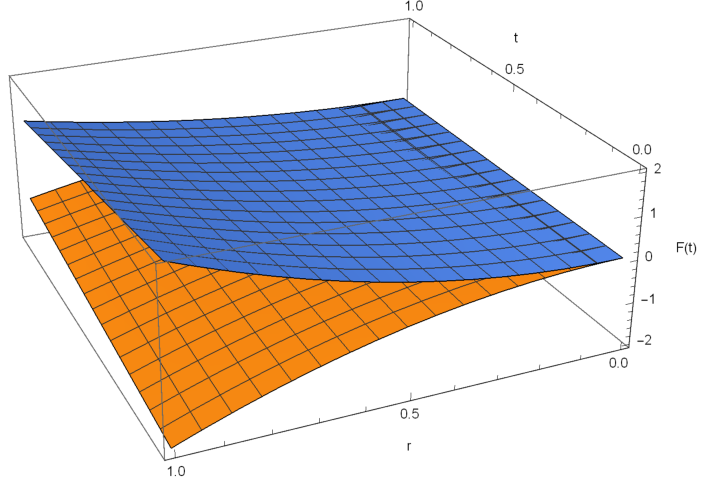}
  \caption{ Graphical plot of fuzzy approximate solution. }\label{f1}
\end{figure}

\begin{figure}[h!]
\centering
 \includegraphics[width=3in]{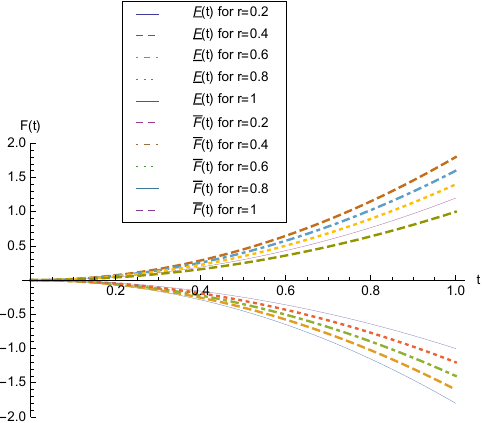}
  \caption{ Fuzzy approximate solution for various $\mu$. }\label{f2}
\end{figure}

\begin{figure}[h!]
\centering
$$\begin{array}{cc}
 \includegraphics[width=3in]{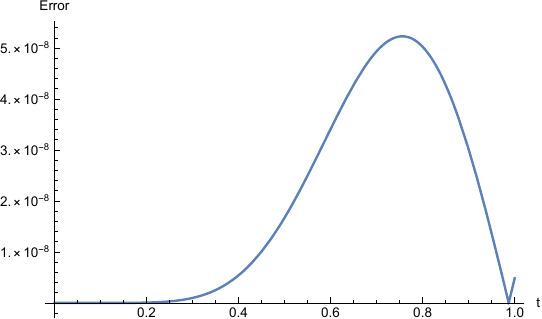}& \includegraphics[width=3.2in]{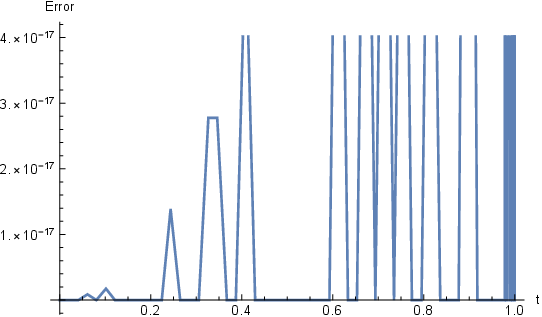}\\
 (a) & (b) \\
 \end{array}$$
  \caption{ Error for (a) $\underline{F}(v), m=5, \mu=0.5$ (b) $\overline{F}(v), m=10, \mu=0.5$. }\label{f3}
\end{figure}

%====================================================================================================================
 \begin{Example}\label{ex1} We consider the problem (\ref{4}) with $K_1(r,v) = 1+v-r, K_2 (r,v) = v-1, m'=2$, $a= \theta_0(v) = 0, \theta_1(v) = \frac{t}{3}$ and $\theta_2(v) = v$ where
 $$
 \begin{array}{l}
 \displaystyle   \underline{Y}  = (-2 + \mu) v^2 - \frac{26}{81} (-2 + \mu) (-1 + v) v^3 -
 \frac{1}{324} (-2 + \mu) v^3 (4 + 3 v),\\
   \\
 \displaystyle   \overline{Y} =  (2 - \mu) v^2 + \frac{26}{81}  (-2 + \mu) (-1 + v) v^3 +
 \frac{1}{324} (-2 + \mu) v^3 (4 + 3 v),
 \end{array}
 $$
and the exact solution $(\underline{F}(v), \overline{F}(v)) = ((\mu - 2) v^2, (2-\mu)v^2)$. In Fig. \ref{f1} the graph of approximate solutions for $(\underline{F}_{20}(v), \overline{F}_{20}(v))$ is presented. Also the graph of obtained solutions for various $r$ and $m=20$  is demonstrated in Fig. \ref{f2}. Fig \ref{f3} shows the error functions for various iterations and $\mu=0.5$.
 \end{Example}

%\begin{table}[h!]
%\caption{ Applying the iterative method with FPA for \ref{ex1} with $\eta=10^{-4}$.}\label{tb1} \centering\scalebox{0.9}{
%\begin{tabular}{|l|l|l|}
 % \hline
%$n$ & $ \tau_{n+1}(\mu)$ ~~~~~~~~~~~~& $|\tau_{n+1}(\mu) -\tau(\mu)|$  ~~~~~~~~~~~~\\
 %   \hline
 %   1      &     1.09445846080780029297    &        0.00037908554077148438\\
 %   2      &     1.09485423564910888672    &        0.00001668930053710938\\
 %   \hline
%  \end{tabular}}
%\end{table}

\begin{figure}[h!]
\centering
$$\begin{array}{cc}
 \includegraphics[width=3in]{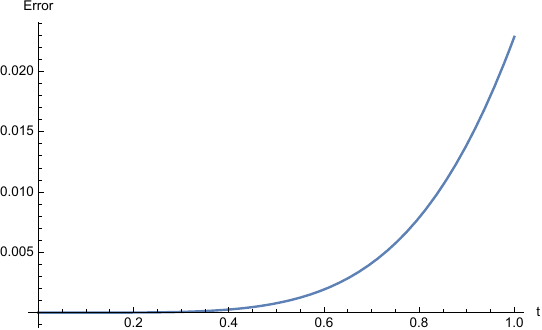}& \includegraphics[width=3.2in]{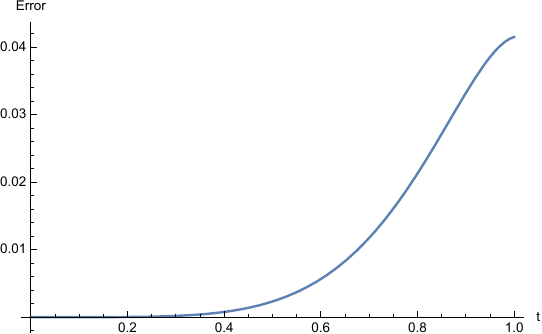}\\
 (a) & (b) \\
 \end{array}$$
  \caption{ Error for (a) $\underline{F}(v), m=20, \mu=0.5$ (b) $\overline{F}(v), m=20, \mu=0.5$. }\label{f5}
\end{figure}

%====================================================================================================================
 \begin{Example}\label{ex2} We have  $K_1(r,v) = v, K_2 (r,v) = v-1, K_3 (r,v) = r-v , m'=3$, $z_1= \theta_0(v) = 0, \theta_1(v) = \frac{t}{8}, \theta_2(v) = \frac{2t}{8}$ and $\theta_3(v) = v$, with nonlinear term $G(Z(r)) = F^3(r)$ where
 $$
 \begin{array}{l}
 \displaystyle   \underline{Y}  = (1 - \mu) v^3 +\frac{ (1023 (-1 + \mu)^3 (-1 + v) v^{10})}{10737418240} - \frac{(
 1073733109 (-1 + \mu)^3 v^{11})}{118111600640},\\
   \\
 \displaystyle   \overline{Y} =  (1 + \mu)v^3 - \frac{(1023 (1 + \mu)^3 (-1 + v) v^{10})}{10737418240} + \frac{(
 1073733109 (1 + \mu)^3 v^{11})}{118111600640},
 \end{array}
 $$
and the exact solution $(\underline{F}(v), \overline{F}(v)) = ((1 - \mu) v^3, (1 + \mu) v^3)$. The graph of obtained solutions for various $\mu$ and $m=10$ is demonstrated in Fig. \ref{f2}. Also the error function for $m=20$ and $\mu=0.5$ can be found in Fig. \ref{f5}.
 \end{Example}

\begin{figure}[h!]
\centering
 \includegraphics[width=3in]{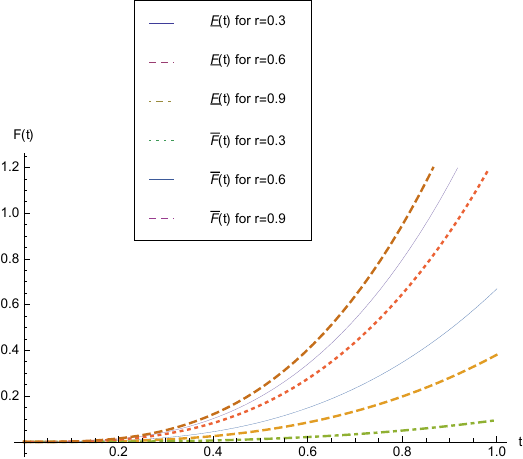}
  \caption{ Fuzzy approximate solution for various $\mu$. }\label{f4}
\end{figure}

~\\
~\\
~\\
~\\
%%%======================================================================================================================
\section{Conclusion}
In this work, the fuzzy Volterra integral equation of the second kind with piecewise kernel was studied. We applied the successive approximation scheme. The existence of an unique solution with the error bound and also the error estimation theorems were discussed. Some examples  have been discussed. Plotting the graphs of fuzzy approximate solutions for various $\mu$ and error functions we showed the accuracy of the method. As the future works, we will combine the method with the CESTAC-CADNA strategy to find the numerical optimality results and optimal distance.
%%%======================================================================================================================

\section*{Funding}
The research was supported by RSF (Project No. 22-29-01619).

%%%======================================================================================================================
%%%======================================================================================================================
%%%======================================================================================================================

% ------------------------------------------------------------------------


\begin{thebibliography}{99}



\bibitem{29a} D. Sidorov,  Integral Dynamical Models: Singularities, Signals and Control. Wold Scientifc. Singapore. 2014.

\bibitem{29aa} A. Apartsyn,  Nonclassical Linear Volterra Equations of the First Kind. Inverse and ill-posed problems series, Vol. 39. Utrecht, Boston: VSP. 2003, 168 p.

\bibitem{m8} B. Bede, S.G. Gal, Quadrature rules for integrals of fuzzy-number-valued functions, Fuzzy Sets and Systems 145 (2004) 359–380.
\bibitem{m13} D. Dubois, H. Prade, Towards fuzzy differential caculus, Fuzzy Sets and Systems 8 (1982) 1–7 (105–116, 225–233).
\bibitem{m18} M. Friedman, M. Ma, A. Kandel, Solutions to fuzzy integral equations with arbitrary kernels, International Journal of Approximate Reasoning 20 (1999) 249–262.
\bibitem{m22} R. Goetschel, W. Voxman, Elementary fuzzy calculus, Fuzzy Sets and Systems 18 (1986) 31–43.
\bibitem{m23} O. Kaleva, Fuzzy differential equations, Fuzzy Sets and Systems 24 (1987) 301–317.
\bibitem{m25} S. Nanda, On integration of fuzzy mappings, Fuzzy Sets and Systems 32 (1989) 95–101.
\bibitem{m30} J.Y. Park, J.U. Jeong, A note on fuzzy integral equations, Fuzzy Sets Systems 108 (1999) 193–200.
\bibitem{m31} J.Y. Park, H.K. Han, Existence and uniqueness theorem for a solution of fuzzy Volterra integral equations, Fuzzy Sets Systems 105 (1999) 481–488.


%===========================================================================

\bibitem{1} S. Ziari, S. Abbasbandy, Open fuzzy quadrature rule for nonlinear fuzzy integral equations with error approximation, Mathematical Researches 7(4) (1400) 781-796

\bibitem{2} S. Farzaneh Javan, S. Abbasbandy, M. A. Fariborzi Araghi, Reproducing Kernel Hilbert space method for solving fuzzy integral equations of the second kind, J. New Researches in Mathematics 8(36) (2022) 29-42.


\bibitem{3}  SS Asari, M Amirfakhrian, S Chakraverty, Application of radial basis functions in solving fuzzy integral equations,
Neural Computing and Applications 31 (10), 6373-6381 2019

\bibitem{4} M Amirfakhrian, K Shakibi, R Rodríguez López, Fuzzy quasi-interpolation solution for Fredholm fuzzy integral equations of second kind,
Soft Computing 21 (15), 4323-4333  2017

\bibitem{5} T Allahviranloo, S Salahshour, Advances in Fuzzy Integral and Differential Equations, Springer 2022

\bibitem{6} T. Allahviranloo, R. Saneifard, R. Saneifard, F. Kiani, S. Noeiaghdam, V. Govindan, The Best Approximation of Generalized Fuzzy Numbers Based on Scaled Metric, Journal of Mathematics, Volume 2022, Article ID 1414415, 8 pages. https://doi.org/10.1155/2022/1414415


\bibitem{13} S. Noeiaghdam, M.A. Fariborzi Araghi, S. Abbasbandy, Valid implementation of Sinc-collocation method to solve the fuzzy Fredholm integral equation, Journal of Computational and Applied Mathematics, 370 (2020) 112632. https://doi.org/10.1016/j.cam.2019.112632

\bibitem{16} M. A. Fariborzi Araghi, S. Noeiaghdam, Homotopy analysis transform method for solving generalized Abel's fuzzy integral equations of the first kind, 4th Iranian Joint Congress on Fuzzy and Intelligent Systems, CFIS 2015, 2016, 7391645. https://doi.org/10.1109/CFIS.2015.7391645

\bibitem{17}  M.A. Fariborzi Araghi, S. Noeiaghdam (2022) Finding Optimal Results in the Homotopy Analysis Method to Solve Fuzzy Integral Equations. In: Allahviranloo T., Salahshour S. (eds) Advances in Fuzzy Integral and Differential Equations. Studies in Fuzziness and Soft Computing, vol 412. Springer, Cham.  \verb"https://doi.org/10.1007/978-3-030-73711-5_7"

\bibitem{18}  Chapter: S. Noeiaghdam, M. A. Fariborzi Araghi, (2021) Application of the CESTAC Method to Find the Optimal Iteration of the Homotopy Analysis Method for Solving Fuzzy Integral Equations. In: Allahviranloo T., Salahshour S., Arica N. (eds) Progress in Intelligent Decision Science. IDS 2020. Advances in Intelligent Systems and Computing, vol 1301. Springer, Cham.   \verb"https://doi.org/10.1007/978-3-030-66501-2_49"

%================================================================================================================================

\bibitem{25}  M. V. Falaleev, N. A. Sidorov, D. N. Sidorov, “Generalized solutions of Volterra integral equations of the first kind”, Lobachevskii J. Math., 20 (2005), 47–57

\bibitem{26} N. A. Sidorov, M. V. Falaleev, D. N. Sidorov, “Generalized solutions of Volterra integral equations of the first kind”, Bulletin of the Malaysian Mathematical Sciences Society, 28:2 (2006), 101–109

\bibitem{27} I. R. Muftahov, D. N. Sidorov, “Solvability and numerical solutions of systems of nonlinear Volterra integral equations of the first kind with piecewise continuous kernels”, Vestn. YuUrGU. Ser. Matem. modelirovanie i programmirovanie, 9:1 (2016), 130–136

\bibitem{28} N. A. Sidorov, D. N. Sidorov, “On the Solvability of a Class of Volterra Operator Equations of the First Kind with Piecewise Continuous Kernels”, Math. Notes, 96:5 (2014), 811–826


\bibitem{29} D. N. Sidorov, “Solvability of systems of integral Volterra equations of the first kind with piecewise continuous kernels”, Russian Math. (Iz. VUZ), 57:1 (2013), 54–63

%================================================================================================================================

\bibitem{31} N.A. Sidorov, D.N. Sidorov, A.V. Krasnik, “Solution of Volterra operator-integral equations in the nonregular case by the successive approximation method”, Differential Equations, 46:6 (2010), 882–891

\bibitem{32} Muftahov I., Tynda A., Sidorov D., “Numeric solution of Volterra integral equations of the first kind with discontinuous kernels”, Journal of Computational and Applied Mathematics, 2017, 119–128


\bibitem{19} A. Tynda, S. Noeiaghdam, D. Sidorov, Polynomial Spline Collocation Method for Solving Weakly Regular Volterra Integral Equations of the First Kind. The Bulletin of Irkutsk State University. Series Mathematics, 2022, vol. 39, pp. 62–79. https://doi.org/10.26516/1997-7670.2022.39.62

\bibitem{20} S. Noeiaghdam, S. Micula, A Novel Method for Solving Second Kind Volterra Integral Equations with Discontinuous Kernel. Mathematics 2021, 9, 2172. https://doi.org/10.3390/math9172172

\bibitem{21} S. Noeiaghdam, D. Sidorov, Integral equations: Theories, Approximations and Applications, Symmetry 2021, 13, 1402. https://doi.org/10.3390/sym13081402

\bibitem{22} S. Noeiaghdam, D. Sidorov, A. M. Wazwaz, N. Sidorov, V. Sizikov, The numerical validation of the Adomian decomposition method for solving Volterra integral equation with discontinuous kernel using the CESTAC method, Mathematics, 2021, 9(3), 1–15, 260. https://doi.org/10.3390/math9030260

\bibitem{23} S. Noeiaghdam, A. Dreglea, J. H. He, Z. Avazzadeh, M. Suleman, M. A. Fariborzi Araghi, D. Sidorov, N. Sidorov, Error estimation of the homotopy perturbation method to solve second kind Volterra integral equations with piecewise smooth kernels: Application of the CADNA library, Symmetry, 2020, 12(10), 1–16, 1730. https://doi.org/10.3390/sym12101730

\bibitem{24} S. Noeiaghdam, D. Sidorov, V. Sizikov, N. Sidorov, Control of accuracy on Taylor-collocation method to solve the weakly regular Volterra integral equations of the first kind by using the CESTAC method, Applied and Computational Mathematics an International Journal, 19 (1) (2020) 81-105.






\end{thebibliography}
\end{document}